\documentclass{article}
\usepackage{amsfonts}

\usepackage{graphicx}
\usepackage{amsmath}
\usepackage[symbol]{footmisc}
\usepackage[french]{babel}

\newtheorem{theorem}{Theorem}[section]

\newtheorem{corollary}[theorem]{Corollary}

\newtheorem{definition}[theorem]{Definition}

\newtheorem{lemma}[theorem]{Lemma}
\newtheorem{notation}[theorem]{Notation}

\newtheorem{proposition}[theorem]{Proposition}
\newtheorem{remark}[theorem]{Remark}

\newenvironment{proof}[1][Proof]{\textbf{#1.} }{\ \rule{0.5em}{0.5em}}
\input{tcilatex}

\begin{document}

\begin{center}
{\huge A q-analogue of certain symmetric functions and one of its
specializations}

\bigskip

{\large Vincent Brugidou} 
\footnotetext{\textit{E-mail address:}vincent.brugidou@univ-lille.fr}
\end{center}

\bigskip

\begin{center}
\textit{Universit\'{e} de Lille, 59655 Villeneuve d'Ascq cedex, France}
\end{center}

\bigskip

\textbf{Abstract: }Let the symmetric functions be defined for the pair of
integers $\left( n,r\right) $, $n\geq r\geq 1$, by $p_{n}^{\left( r\right)
}=\sum m_{\lambda }$ where $m_{\lambda }$ are the monomial symmetric
functions, the sum being over the partitions $\lambda $ of the integer $n$
with length $r$. We introduce by a generating function, a $q$-analogue of $%
p_{n}^{\left( r\right) }$ and give some of its properties. This $q$-analogue
is related to its the classical form using the $q$-Stirling numbers. We also
start with the same \ procedure the study of a $p,q$-analogue of $%
p_{n}^{\left( r\right) }$.

By specialization of this $q$-analogue in the series $\sum\nolimits_{n=0}^{%
\infty }q^{\binom{n}{2}}t^{n}/n!$, we recover in a purely formal way$\ $a
class of polynomials $J_{n}^{\left( r\right) }$ historically introduced as
combinatorial enumerators, in particular of tree inversions. This also
results in a new linear recurrence for those polynomials whose triangular
table can be constructed, row by row, from the initial conditions $%
J_{r}^{\left( r\right) }=1$. The form of this recurrence is also given for
the reciprocal polynomials of \ $J_{n}^{\left( r\right) }$, known to be the
sum enumerators of parking functions. Explicit formulas for $J_{n}^{\left(
r\right) }$ and their reciprocals are deduced, leading inversely to new
representations of these polynomials as forest statistics.

\bigskip

\textit{keywords: }Symmetric functions, $q$-analogue, $q$-Stirling numbers,
tree inversions, parking functions.

.

\section{\protect\bigskip \textbf{Introduction}}

This paper is the first in a series whose object is a $q$-analogue defined
in a fairly natural way of certain symmetric functions namely the functions
defined for each pair of integers $\left( n,r\right) $ $\ $such that $n\geq
r\geq 1$ by 
\begin{equation}
p_{n}^{(r)}=\sum\limits_{\left| \lambda \right| =n,\;l\left( \lambda \right)
=r}m_{\lambda }\text{.}  \tag{1.1}
\end{equation}
In $\left( 1.1\right) $ the $m_{\lambda \text{ }}$ are the monomial
symmetric functions, the sum being over the integer partitions $\lambda $ of 
$n$, of length $l(\lambda )=r$. These functions $p_{n}^{(r)}$ are introduced
with this notation in [12, Example 19 p. 33], whose notations we will follow
faithfully. The $q$-analogues thus defined, that we will note $\left[
p_{n}^{(r)}\right] $, have attractive properties. The article presents the
definition and some properties of these $q$-analogues. We also present some
applications by specialization in the following formal series sometimes
called $q$-deformation of the exponential series : 
\begin{equation}
E_{xq}(t)=\sum\limits_{n=0}^{\infty }q^{\binom{n}{2}}\dfrac{t^{n}}{n!}\text{.%
}  \tag{1.2}
\end{equation}

This allows new insight into polynomials known as inversion enumerators in
rooted forests, or in a reciprocal way as the sum enumerators of parking
functions (see Yan's summary article [21] on these notions). New identities
are deduced for these polynomials, from which we extract some combinatorial
consequences.

The organization of the article is as follows. In Section 2 we set the
notations and recall the necessary prerequisites for symmetric functions,
integer partitions and $q$-calculus.

Sections 3 and 4 are devoted to the general study of the $\left[ p_{n}^{(r)}%
\right] $. In Section 3, starting from a generating function of $%
p_{n}^{\left( r\right) }$ we provide the definition of $\left[ p_{n}^{(r)}%
\right] $\ along with its initial properties, including the particular case $%
r=1$. Section 4 contains Theorem 4.1, which is one of the main results of
the article. It establishes a connection between the $q$-analogue and its
classical form using the $q$-Stirling numbers of the second kind $S_{q}\left[
j,r\right] $, through the following equation (Equation $\left( 4.4\right) $):

\begin{equation*}
\left[ p_{n}^{(r)}\right] _{q}=\sum\limits_{j=r}^{n}\left( 1-q\right)
^{j-r}S_{q}\left[ j,r\right] \,p_{n}^{(j)}\text{,}
\end{equation*}

We also give the matrix form and the inverse equations of this theorem
through $q$-Stirling numbers of the first kind. Section 5 begins the
generalization to $p,q$-analogue.

The remaining Sections 6 to 9 are devoted to the applications of $\left[
p_{n}^{(r)}\right] $, obtained through specialization in $E_{xq(t)}$. In
Section 6, Theorem 4.1 leads us to polynomials previously introduced in $%
\left[ 16\right] $ and [20], but in a different way from that used by those
authors. More precisely, it concerns the class of polynomials defined for $%
n\geq r\geq 1$\ that we will denote by $J_{n}^{\left( r\right) }$, and which
, according to the notations in [21], correspond to the polynomials $%
I_{n-r}^{(r,1)}$.

A new result that we obtained through this approach is the following linear
recurrence for these polynomials (Equation $\left( 6.5\right) $):

\begin{equation*}
\text{For }n-1\geq r\geq 1\text{, \ \ \ \ \ \ }J_{n}^{\left( r\right)
}(q)=\sum\limits_{j=1}^{n-r}\left[ r\right] _{q}^{j}\;q^{\binom{j}{2}}\binom{%
n-r}{j}J_{n-r}^{\left( j\right) }(q).
\end{equation*}
This recurrence allows all the polynomials to be determined from the initial
conditions $J_{r}^{\left( r\right) }=1$. As an example, the first rows and
columns of the table of $J_{n}^{\left( r\right) }$ are provided.

In Section 7, we study the reciprocal polynomials of $J_{n}^{\left( r\right)
}$, denoted by $\overline{J_{n}^{\left( r\right) }}$, and which are sum
enumerators of parking functions. Recurrence $\left( 6.5\right) $\
translates into a linear recurrence between $\overline{J_{n}^{\left(
r\right) }}$, that we compare with another linear recurrence relation coming
from the application of Goncarev polynomials to parking functions [11].

In Section 8, by applying Recurrence $\left( 6.5\right) $, we derive
explicit formulas for the $J_{n}^{\left( r\right) }$ polynomials presented
in Theorem 8.2, which represent additional new results of our work. In
Section 9, we give combinatorial interpretations of these explicit formulas
by introducing new statistics on forests whose enumerator polynomials are $%
J_{n}^{\left( r\right) }$ or $\overline{J_{n}^{\left( r\right) }}$. \ In
particular, this gives a $q$-refinement of the enumeration of functionnal
digraph given in [14, page 19].

Complementary articles will follow, in particular $\left[ 2\right] $,
presenting other properties or applications of this $q$-analogue.

\section{Preliminaries}

\QTP{Body Math}
Let us first recall the definitions relating to integer partitions and
symmetric functions, refering for more details to [12, Chap. 1] or [17,
Chap. 7]. If $\lambda =\left( \lambda _{1},\lambda _{2},...,\lambda
_{r}\right) $ is a partition of the integer $n$, we set $\left| \lambda
\right| =\lambda _{1}+\lambda _{2}+...+\lambda _{r}=n,\;l(\lambda )=r$, and$%
\;$

\QTP{Body Math}
\begin{equation}
n(\lambda )=\sum\limits_{i\geq 1}\left( i-1\right) \lambda
_{i}=\sum\limits_{i\geq 1}\binom{\lambda _{i}^{\prime }}{2}\text{,} 
\tag{2.1}
\end{equation}
\ where $\lambda ^{\prime }=\left( \lambda _{1}^{\prime },\lambda
_{2}^{\prime },..\right) $ is the conjugate of the partition $\lambda $.\ 

\QTP{Body Math}
If $A$ is a commutative ring, $\Lambda _{A}$ is the set of symmetric
functions in the variables $X=\left( x_{i}\right) _{i\geq 1}\;$with
coefficients in $A$. For our purposes $A$ will be the field of rational
functions $\mathbb{Q}\left( q\right) ,$ with eventually some indeterminates
adjoined (in particular $p$, at Section 5). For $\lambda $ describing the
set of partitions, $\left( m_{\lambda }\right) ,\left( e_{\lambda }\right)
,\left( h_{\lambda }\right) $ and $\left( p_{\lambda }\right) $, are the
four classical bases of $\Lambda _{A}$. Agreeing that $e_{0}=h_{0}=1\;$we
recall that

\QTP{Body Math}
\begin{equation}
E(t)=\sum\limits_{n\geq 0}e_{n}t^{n}=\prod\limits_{i\geq
1}(1+x_{i}t)\;;\;\;H(t)=\sum\limits_{n\geq 0}h_{n}t^{n}=\left( E(-t)\right)
^{-1}=\dfrac{1}{\prod\limits_{i\geq 1}\left( 1-x_{i}t\right) }\text{.} 
\tag{2.2}
\end{equation}

\QTP{Body Math}
For the $q$-analogue we use, with a few exceptions, the notations of [9], to
which we refer for more details. The index $q$ can be omitted if there \ is
no ambiguity. For ($n,k)\in \mathbf{N}^{2}$, we have

\QTP{Body Math}
\begin{equation*}
\left[ n\right] _{q}=1+q+q^{2}+...+q^{n-1}\text{ for }n\neq 0\;\text{and }%
\left[ 0\right] _{q}=0\text{,}
\end{equation*}

\QTP{Body Math}
\begin{equation*}
\left[ n\right] _{q}!=\left[ 1\right] \left[ 2\right] ...\left[ n\right] 
\text{ for }n\neq 0\;\text{and }\left[ 0\right] _{q}!=1\text{,}
\end{equation*}

\QTP{Body Math}
\begin{equation*}
\QATOPD[ ] {n}{k}_{q}=\dfrac{\left[ n\right] !}{\left[ k\right] !\left[ n-k%
\right] !}\text{ for }n\geq k\geq 0\text{ \ \ and }\QATOPD[ ] {n}{k}_{q}=0\;%
\text{otherwise.}
\end{equation*}

\QTP{Body Math}
$\bigskip $The $q$-derivation of a formal series $F(t)$ is given by 
\begin{equation*}
D_{q}F(t)=\dfrac{F(qt)-F(t)}{(q-1)t}\text{ \ \ and for }r\geq
1\;\;D_{q}^{r}=D_{q}\left( D_{q}^{r-1}F(t)\right) \text{ \ with \ }%
D_{q}^{0}F(t)=F\text{.}
\end{equation*}

\QTP{Body Math}
Specifically 
\begin{equation}
D_{q}^{r}t^{n}=\left[ r\right] !\QATOPD[ ] {n}{r}t^{n-r}\text{ for }n\geq r\
\ \text{and}\ \ D_{q}t^{0}=0\text{.}  \tag{2.3}
\end{equation}

\QTP{Body Math}
$\bigskip $

\QTP{Body Math}
We also use as an index instead of $q$, a rational fraction of $q$ with
coefficients in $\mathbb{Q}$. It is necessary to understand in this case
what we will give as a definition.

\begin{definition}
Let $\psi (q)$ $\in \mathbb{Q}\left( q\right) $ such that $\psi (1)=1$. We
set for all $n\in \mathbb{N}$%
\begin{equation*}
\left[ n\right] _{\psi (q)}=1+\psi (q)+\left( \psi (q)\right)
^{2}+...+\left( \psi (q)\right) ^{n-1}\ \ \text{if}\ n\geq 1\ \ \text{and}\
\ \ \left[ 0\right] _{\psi (q)}=0\text{,}
\end{equation*}

then are defined as above $\left[ n\right] _{\psi (q)}!$, $\QATOPD[ ] {n}{k}%
_{\psi (q)}$, etc... . In the same way we set 
\begin{equation*}
D_{\psi (q)}F(t)=\dfrac{F(\psi (q)t)-F(t)}{(\psi (q)-1)t}\text{.}
\end{equation*}
\end{definition}

\QTP{Body Math}
Finally, some usual notations used in the article are recalled. $D^{r}F(t)$
is the usual derivation of order $r$ with respect to $t$ of the formal
series $F(t)$, $\delta _{i}^{j}=\delta _{i,j}$ is the Kronecker symbol which
is 1 if $i=j$ and 0 otherwise. If $E$ is a finite set, $\left| E\right| $
denotes its cardinality and $\mathbf{2}^{E}$ the set of parts of $E$. $%
\mathbb{N},\mathbb{Z}$ denote respectively the sets of natural and relative
integers and $\mathbb{Q}$\ that of rational numbers.$\ $If $A$ is a
commutative ring and $t$ an indeterminate, $A\left[ t\right] $ and $A\left[ %
\left[ t\right] \right] $) are the set of polynomials and the set of formal
series respectively, in $t$ with coefficients in $A$. For $n\in \mathbb{N}%
^{\ast }=\mathbb{N}-\left\{ 0\right\} $ we denote $\mathbf{n=}\left\{
1,2,...n\right\} $.

\section{Definition of $\left[ p_{n}^{(r)}\right] $}

\bigskip

\subsection{\protect\bigskip Some properties of $p_{n}^{\left( r\right) }$}

Let us give some properties of $p_{n}^{(r)}$ defined by $\left( 1.1\right) $%
, which are given in [12, example 19, p.33], or can easily be deduced from
this exemple. By agreeing to set $p_{n}^{(0)}=\delta _{n}^{0}$, the second
equation from [12, example 19 p.33] gives, since $H\left( t\right) =\left(
E\left( -t\right) \right) ^{-1}$, \ the following generating function of $%
p_{n}^{\left( r\right) }$ for all $r\in \mathbf{N}$: 
\begin{equation}
\text{\ }\sum\limits_{n\geq r}p_{n}^{(r)}t^{n-r}=\dfrac{\left( -1\right) ^{r}%
}{r!}\dfrac{D^{r}\left( E\left( -t\right) \right) }{E\left( -t\right) }\text{%
,}  \tag{3.1}
\end{equation}
where $D^{r}\left( E\left( -t\right) \right) $ is the $r-$th derivative of $%
E\left( -t\right) $ with respect to $t$.

By making the change of variable $u=-t$, $\left( 3.1\right) $ becomes
equivalent to:

\begin{equation}
\sum\limits_{n\geq r}p_{n}^{(r)}\left( -t\right) ^{n-r}=\dfrac{1}{r!}\dfrac{%
D^{r}E\left( t\right) }{E(t)}\text{.}  \tag{3.2}
\end{equation}
Since $D^{r}\left( t^{n}\right) =r!\binom{n}{r}t^{n-r}$, we have by
linearity: 
\begin{equation}
D^{r}E(t)=r!\sum\limits_{n\geq r}\binom{n}{r}e_{n}t^{n-r}\text{.}  \tag{3.3}
\end{equation}
From $\left( 3.2\right) $ and $\left( 3.3\right) $, it follows that: 
\begin{equation*}
\sum\limits_{n\geq 0}e_{n}t^{n}\sum\limits_{n\geq
r}p_{n}^{(r)}(-t)^{n-r}=\sum\limits_{n\geq r}\binom{n}{r}e_{n}t^{n-r}\text{.}
\end{equation*}
\qquad By equating the coefficients, we obtain a system of linear equations
for the unknowns $p_{n}^{\left( r\right) }$: 
\begin{equation*}
\text{for }n\geq r\text{,}\;\;\;\sum\limits_{k=r}^{n}\left( -1\right)
^{k-r}e_{n-k}p_{k}^{(r)}=\binom{n}{r}e_{n}\text{,}
\end{equation*}
for which the following solutions are given by Cramer's rule: 
\begin{equation}
p_{n}^{(r)}=\left| 
\begin{array}{cccccccc}
\binom{r}{r}e_{r} & 1 & 0 & 0 & . & . & . & 0 \\ 
\binom{r+1}{r}e_{r+1} & e_{1} & 1 & 0 & . & . & . & 0 \\ 
\binom{r+2}{r}e_{r+2} & e_{2} & e_{1} & 1 & 0 & . & . & 0 \\ 
. & . &  & . & . & . & . & . \\ 
. & . &  &  & . & . & . & . \\ 
. & . &  &  &  & . & . & 0 \\ 
. & . &  &  &  &  & . & 1 \\ 
\binom{n}{r}e_{n} & e_{n-r} & e_{n-r-1} & . & . & . & e_{2} & e_{1_{{}}}
\end{array}
\right| \text{.}  \tag{3.4}
\end{equation}

\subsection{\protect\bigskip Definition of $\left[ p_{n}^{\left( r\right) }%
\right] $}

To define $\left[ p_{n}^{(r)}\right] _{q}$, we will use the following $q$%
-analogue of $\left( 3.1\right) $.

\begin{definition}
For any pair of integers $\left( n,r\right) $\ such that $n\geq r\geq 0$, we
set: 
\begin{equation}
\sum\limits_{n\geq r}\left[ p_{n}^{(r)}\right] _{q}t^{n-r}=\dfrac{\left(
-1\right) ^{r}}{\left[ r\right] _{q}!}\dfrac{D_{q}^{r}\left( E(-t)\right) }{%
E(-t)}\text{,}  \tag{3.1q}
\end{equation}
which implies in particular $\left[ p_{n}^{(0)}\right] _{q}=\delta _{n}^{0}$.
\end{definition}

The subscript $q$ will sometimes be omitted in the following if there is no
ambiguity. It is easy to verify that the calculations performed in
Subsection 3.1 can be transposed to this $q$-analogue of $p_{n}^{\left(
r\right) }$. In particular $\left( 3.1q\right) $ is equivalent to: 
\begin{equation}
\sum\limits_{n\geq r}\left[ p_{n}^{(r)}\right] \left( -t\right) ^{n-r}=%
\dfrac{1}{\left[ r\right] !}\dfrac{D_{q}^{r}E\left( t\right) }{E(t)}\text{,}
\tag{3.2q}
\end{equation}

and 
\begin{equation}
D_{q}^{r}E(t)=\left[ r\right] !\sum\limits_{n\geq r}\QATOPD[ ] {n}{r}%
e_{n}t^{n-r}\text{.}  \tag{3.3q}
\end{equation}
\bigskip Thus, we have the system of linear equations:

\begin{equation*}
\text{For }n\geq r\;\;\;\sum\limits_{k=r}^{n}\left( -1\right) ^{k-r}e_{n-k}%
\left[ p_{k}^{(r)}\right] =\QATOPD[ ] {n}{r}e_{n}\text{,}
\end{equation*}
whose solutions is the $q$-analogue of $\left( 3.4\right) $, given by the
following proposition:

\begin{proposition}
For $n\geq r\geq 1$, we have 
\begin{equation}
\left[ p_{n}^{(r)}\right] =\left| 
\begin{array}{cccccccc}
\QATOPD[ ] {r}{r}e_{r} & 1 & 0 & 0 & . & . & . & 0 \\ 
\QATOPD[ ] {r+1}{r}e_{r+1} & e_{1} & 1 & 0 & . & . & . & 0 \\ 
\QATOPD[ ] {r+2}{r}e_{r+2} & e_{2} & e_{1} & 1 & 0 & . & . & 0 \\ 
. & . &  & . & . & . & . & . \\ 
. & . &  &  & . & . & . & . \\ 
. & . &  &  &  & . & . & 0 \\ 
. & . &  &  &  &  & . & 1 \\ 
\QATOPD[ ] {n}{r}e_{n} & e_{n-r} & e_{n-r-1} & . & . & . & e_{2} & e_{1_{{}}}
\end{array}
\right| \text{.}  \tag{3.4q}
\end{equation}
\end{proposition}

This determinant can be taken as an alternative definition of $\left[
p_{n}^{(r)}\right] $.

\QTP{Body Math}
$\bigskip $

\QTP{Body Math}
\textbf{Particular case} $r=1$

Equation $\left( 1.1\right) $ gives in this case $p_{n}^{\left( 1\right)
}=m_{n}=p_{n}$ thus $\left[ p_{n}^{\left( 1\right) }\right] $ can be also
denoted $\left[ p_{n}\right] $. So, we have from $\left( 3.2q\right) $

\begin{equation*}
\sum\limits_{n\geq 1}\left[ p_{n}\right] _{q}\left( -t\right) ^{n-1}=\dfrac{%
D_{q}E(t)}{E(t)}\text{,}
\end{equation*}
and the system of linear equations 
\begin{equation}
\text{ For }n\geq 1\;\;\;\sum\nolimits_{k=1}^{n}\left( -1\right)
^{k-1}e_{n-k}\left[ p_{k}\right] =\left[ n\right] e_{n}\text{,}  \tag{3.5}
\end{equation}
whose solution are

\begin{equation}
\left[ p_{n}\right] =\left| 
\begin{array}{ccccccc}
\left[ 1\right] e_{1} & 1 & 0 & 0 & . & . & 0 \\ 
\left[ 2\right] e_{2} & e_{1} & 1 & 0 & . & . & 0 \\ 
. & . & . & . & . & . & . \\ 
. & . & . & . & . & . & . \\ 
. & . & . &  & . & . & 0 \\ 
. & . & . &  &  & . & 1 \\ 
\left[ n\right] e_{n} & e_{n-1} & e_{n-2} & . & . & e_{2} & e_{1}
\end{array}
\right| \text{.}  \tag{3.6}
\end{equation}
By inverting $\left( 3.5\right) $, we obtain the system of linear equations
in the unknwons $e_{n}$ 
\begin{equation*}
\text{ for }n\geq 1\text{,}\;\;\;\sum\nolimits_{k=1}^{n-1}\left( -1\right)
^{k-1}\left[ p_{n-k}\right] e_{k}+\left[ n\right] e_{n}=\left[ p_{n}\right] 
\text{,}
\end{equation*}
which gives by Cramer's rule:

\begin{equation}
\left[ n\right] !e_{n}=\left| 
\begin{array}{ccccccc}
\left[ p_{1}\right] & \left[ 1\right] & 0 & 0 & . & . & 0 \\ 
\left[ p_{2}\right] & \left[ p_{1}\right] & \left[ 2\right] & 0 & . & . & 0
\\ 
. & . & . & . & . & . & . \\ 
. & . & . & . & . & . & . \\ 
. & . &  & . & . & . & 0 \\ 
. & . &  &  &  & . & \left[ n-1\right] \\ 
\left[ p_{n}\right] & \left[ p_{n-1}\right] & \left[ p_{n-2}\right] & . & .
& \left[ p_{2}\right] & \left[ p_{1}\right]
\end{array}
\right| \text{.}  \tag{3.7}
\end{equation}

We notice that $\left( 3.6\right) $ and $\left( 3.7\right) $\ are $q$%
-analogues of the classical case (see the first two equations of [12,
Example 8 p. 28]. In $\left[ 2\right] $, the study of $\left[ p_{n}\right] $
has been developed, and more generally, of: 
\begin{equation}
\left[ p_{\lambda }\right] =\left[ p_{\lambda _{1}}\right] \left[ p_{\lambda
_{2}}\right] ...\left[ p_{\lambda _{r}}\right] \text{,}  \tag{3.8}
\end{equation}
defined for any integer partition $\lambda =\left( \lambda _{1},\lambda
_{2},...,\lambda _{r}\right) $.

\section{\protect\bigskip Relations between $\left[ p_{n}^{(r)}\right] $ and 
$p_{n}^{(r)}$}

To prove Theorem 4.1, some reminders about the $q$-analogues of the Stirling
numbers must be made. The classical Stirling numbers are well known (see,
for example, [5, Chap. 5]). For their $q$-analogues which are still the
subject of research, references are made to [6, 15]. In [3, 4], Carlitz
defined the $q$-Stirling numbers of the second kind, which we denote\ $S_{q}%
\left[ n,k\right] $, by the following identity (with our notations): 
\begin{equation}
\QATOPD[ ] {n}{k}_{q}=\sum\limits_{j=k}^{n}\binom{n}{j}\left( q-1\right)
^{j-k}S_{q}\left[ j,k\right] \text{,}  \tag{4.1}
\end{equation}
(the index $q$ is omitted if there is no ambiguity). By inverted $\left(
4.1\right) $, Carlitz obtained:

\begin{equation}
\left( 1-q\right) ^{n-k}S\left[ n,k\right] =\sum\limits_{l=k}^{n}\left(
-1\right) ^{l-k}\binom{n}{l}\QATOPD[ ] {l}{k}\text{,}  \tag{4.2}
\end{equation}
as well as the recurrence relation 
\begin{equation}
S\left[ n,k\right] =S\left[ n-1,k-1\right] +\left[ k\right] S\left[ n-1,k%
\right] \text{,}  \tag{4.3}
\end{equation}
which makes it possible to find all the values of $S\left[ n,k\right] $ for $%
\left( n,k\right) \in \mathbf{N}^{2}$ by setting $S\left[ n,0\right] =\delta
_{n}^{0}$,$\;$and$\;S\left[ n,k\right] =0$ if $k>n$. We have $S\left[ n,n%
\right] =S\left[ n,1\right] =1$ and\ \ $S_{q=1}\left[ n,k\right] =S(n,k)$,
where $S(n,k)$ is the classical Stirling numbers of the second kind.
Equation $\left( 4.3\right) $ \ is the $q$-analogue of the recurrence
formula for $S(n,k)$. Note that $\left( 4.1\right) $ and $\left( 4.2\right) $
do not have classical correspondents.

\begin{theorem}
\bigskip We have for all $n\geq r\geq 1$: 
\begin{equation}
\left[ p_{n}^{(r)}\right] _{q}=\sum\limits_{j=r}^{n}\left( 1-q\right)
^{j-r}S_{q}\left[ j,r\right] \,p_{n}^{(j)}\text{,}  \tag{4.4}
\end{equation}
where the $S_{q}\left[ n,r\right] $ are the q-Stirling numbers of the second
kind by Carlitz.
\end{theorem}

\begin{proof}
We start from Equation $\left( 3.2q\right) $ which becomes with $E(t)H(-t)=1$
and $m=n-r$ : 
\begin{equation*}
\sum\limits_{m\geq 0}\left[ p_{m+r}^{(r)}\right] \left( -1\right)
^{m}t^{m}=\left( \dfrac{1}{\left[ r\right] !}D_{q}^{r}E(t)\right) H(-t)\text{%
.}
\end{equation*}

By replacing the expression in the parentheses on the right-hand side with $%
\left( 3.3q\right) $ and $H(-t)=\sum\nolimits_{l\geq 0}h_{l}\left( -1\right)
^{l}t^{l}$, we obtain by equating the coefficients: 
\begin{equation}
\left[ p_{m+r}^{(r)}\right] =\QATOPD[ ] {r}{r}e_{r}h_{m}-\QATOPD[ ] {r+1}{r}%
e_{r+1}h_{m-1}+...+\left( -1\right) ^{m-1}\QATOPD[ ] {r+m-1}{r}%
e_{r+m-1}+\left( -1\right) ^{m}\QATOPD[ ] {r+m}{r}e_{r+m}\text{.}  \tag{4.5}
\end{equation}

The definition of complete symmetric functions implies that $%
h_{n}=\sum\nolimits_{r=1}^{n}\sum\nolimits_{\left| \lambda \right|
=n,l\left( \lambda \right) =r}m_{\lambda
}=\sum\nolimits_{l=1}^{n}p_{n}^{\left( l\right) }$. Therefore, in general,
we have: 
\begin{equation}
e_{r}h_{n}=e_{r}\sum\limits_{l=1}^{n}p_{n}^{\left( l\right)
}=\sum\limits_{j=r}^{r+n}\binom{j}{r}p_{r+n}^{\left( j\right) }.  \tag{4.6}
\end{equation}
Let us clarify where the term on the right-hand side of\ $\left( 4.6\right) $
comes from. The product of a monomial of $e_{r}$ and a monomial of $h_{n}$
gives a monomial of degree $r+n$, whose number $j$\ of distinct variables is
between $r$ and $r+n$ (for an infinity of variable $x_{i}$, or at least if
the number of variables is greater than or equal to $r+n$). This is
precisely the form of the monomials in the right-hand side. Moreover, in the
left-hand side of $\left( 4.6\right) $, the number of identical monomials of
degree $r+n$ with $j$ distinct variables corresponds to the number of
monomials of $e_{r}$ where $r$ distinct variables are taken from the $j$
letters, that is to say $\binom{j}{r}$.

By using Equation $\left( 4.6\right) $ for each products $e_{r+k}h_{m-k}\;$%
in $\left( 4.5\right) $, it follows that: 
\begin{equation*}
\left[ p_{m+r}^{(r)}\right] =\sum\limits_{k=0}^{m}\left( -1\right) ^{k}%
\QATOPD[ ] {r+k}{r}\sum\limits_{j=r+k}^{m+r}\binom{j}{r+k}p_{m+r}^{(j)}\text{%
,}
\end{equation*}
then, by inverting the sums: 
\begin{equation*}
\left[ p_{m+r}^{(r)}\right] =\sum\limits_{j=r}^{m+r}p_{m+r}^{\left( j\right)
}\sum\limits_{k=0}^{j-r}\left( -1\right) ^{k}\binom{j}{r+k}\QATOPD[ ] {r+k}{r%
}\text{.}
\end{equation*}

With the change of index $l=r+k$ and $n=m+r$, we find 
\begin{equation}
\left[ p_{n}^{(r)}\right] _{q}=\sum\limits_{j=r}^{n}p_{n}^{(j)}\sum%
\limits_{l=r}^{j}\left( -1\right) ^{l-r}\binom{j}{l}\QATOPD[ ] {l}{r}_{q}, 
\tag{4.7}
\end{equation}
which finally leads to Equation $\left( 4.4\right) $ when taking into
account $\left( 4.2\right) .$
\end{proof}

\bigskip

\textbf{Particular case} $r=1$. We get for all $n\geq 1$, 
\begin{equation}
\left[ p_{n}\right] =\sum\limits_{j=1}^{n}\left( 1-q\right) ^{j-1}p_{n}^{(j)}%
\text{.}  \tag{4.8}
\end{equation}

\bigskip

\textbf{Matrix form and inverse}. Since $S\left[ n,k\right] =0$ for $k<n$, $%
\left( 4.4\right) $ can still be written for $n\geq r\geq 1$ 
\begin{equation}
\left[ p_{n}^{(r)}\right] =\sum\limits_{j=1}^{n}\left( 1-q\right) ^{j-r}S%
\left[ j,r\right] \,p_{n}^{(j)}\text{.}  \tag{4.9}
\end{equation}
Let the row matrices be defined by \ $\left[ P_{n}\right] =\left( \left[
p_{n}^{(r)}\right] \right) _{r=1}^{n}\;$and$\;P_{n}=\left( p_{n}^{\left(
r\right) }\right) _{r=1}^{n}$. System $\left( 4.9\right) $ for $n\geq r\geq
1 $, can be written in matrix form\ \ $\left[ P_{n}\right] =P_{n}A_{n}$,\
where $A_{n}$ is the triangular matrix $\left( A_{i,j}\right) _{i,j=1}^{n}$
with $A_{i,j}=\left( 1-q\right) ^{i-j}S\left[ i,j\right] $. It is easy to
see that we have \ 
\begin{equation}
A_{n}=U_{n}\left[ S_{n}\right] U_{n}^{-1}\text{,}  \tag{4.10}
\end{equation}
\ where $\left[ S_{n}\right] $ is the triangular matrix of the $q$-Stirling
numbers of the second kind$\ \left[ S_{n}\right] =\left( S_{q}\left[ i,j%
\right] \right) _{i,j=1}^{n}$ and $U_{n}$ is the diagonal matrix \ 
\begin{equation*}
U_{n}=\left( 
\begin{array}{cccc}
\left( 1-q\right) ^{0} &  &  &  \\ 
& \left( 1-q\right) ^{1} &  &  \\ 
&  & ... &  \\ 
&  &  & \left( 1-q\right) ^{n-1}
\end{array}
\right) \text{.}
\end{equation*}

The Stirling numbers of the first kind, denoted $s(n,k)$, have as $q$%
-analogue the $q$-Stirling numbers of the first kind introduced in [8],
which we will denote by $s_{q}\left[ n,k\right] $. Let the matrix $\left[
s_{n}\right] $ be defined by 
\begin{equation*}
\;\;\;\;\ \left[ s_{n}\right] =\left( s_{q}\left[ i,j\right] \right)
_{i,j=1}^{n}\text{,}
\end{equation*}
then we know (see for example [6, p. 96]) that $\left[ S_{n}\right] $ and $%
\left[ s_{n}\right] $ are inverses of each other. We deduce with $\left(
4.10\right) $, that $A_{n}$ is inversible and that $A_{n}^{-1}=U_{n}\left[
s_{n}\right] U_{n}^{-1}$and $P_{n}=\left[ P_{n}\right] A_{n}^{-1}\ $with$\
A_{n}^{-1}=\left( B_{i,j}\right) _{i,j=1}^{n}$and\ $B_{i,j}=\left(
1-q\right) ^{i-j}s\left[ i,j\right] $. Hence, the corollary equivalent to
Theorem 4.1:

\begin{corollary}
We have for all $n\geq r\geq 1$, 
\begin{equation*}
p_{n}^{(r)}=\sum\limits_{j=r}^{n}\left( 1-q\right) ^{j-r}s_{q}\left[ j,r%
\right] \left[ p_{n}^{(j)}\right] _{q}\text{,}
\end{equation*}
where the $s_{q}\left[ n,r\right] $ are the q-Stirling numbers of the first
kind.
\end{corollary}

\textbf{Remark. }It is easy to see that $\left( 4.9\right) $\ generalizes
for $n\geq r\geq 0$, to $\left[ p_{n}^{(r)}\right] =\sum\nolimits_{j=0}^{n}%
\left( 1-q\right) ^{j-r}S\left[ j,r\right] \,p_{n}^{(j)}$. One could easily
deduce the generalizations of the matrix and inverse forms of this system,
with the augmented $q$-Stirling matrices $\widehat{\left[ S_{n}\right] }%
=\left( S\left[ i,j\right] \right) _{i,j=0}^{n}$ and$\;\widehat{\left[ s_{n}%
\right] }=\left( s\left[ i,j\right] \right) _{i,j=0}^{n}$.\ \ \ \ 

\section{\protect\bigskip Extension to $p,q$-analogue}

We know that there is a calculus with two parameters, denoted $p,q$-analogue
calculus, whose origin dates back at least to 1991, and which is reduced to
the $q$-analogue when $p=1$. Let us briefly recall the definitions that
generalize those of the $q$-analogue given in Section 2. One can consult for
more details [6, 7]. 
\begin{equation*}
\left[ n\right] _{p,q}=\dfrac{p^{n}-q^{n}}{p-q}%
=p^{n-1}+p^{n-2}q+...+pq^{n-2}+q^{n-1}\text{,}
\end{equation*}

\begin{equation*}
\left[ n\right] _{p,q}!=\left[ 1\right] _{p,q}\left[ 2\right] _{p,q}...\left[
n\right] _{p,q}\;\ \text{and}\;\QATOPD[ ] {n}{k}_{p,q}=\dfrac{\left[ n\right]
_{p,q}!}{\left[ k\right] _{p,q}!\left[ n-k\right] _{p,q}!}\text{,}
\end{equation*}

\begin{equation*}
D_{p,q}F(t)=\dfrac{F(pt)-F(qt)}{\left( p-q\right) t}\text{ \ then }%
D_{p,q}^{r}F(t)=D_{p,q}\left( D_{p,q}^{r-1}F(t)\right) \text{.}
\end{equation*}

Specifically 
\begin{equation*}
D_{p,q}^{r}t^{n}=\left[ r\right] _{p,q}!\QATOPD[ ] {n}{r}_{p,q}t^{n-r}\text{.%
}
\end{equation*}

The definition of the p,q-analogue of the $p_{n}^{(r)}$ follows that of the $%
q$-analogue.

\begin{definition}
We set for $n\geq r\geq 0$, 
\begin{equation*}
\sum\limits_{n\geq r}\left[ p_{n}^{(r)}\right] _{p,q}t^{n-r}=\dfrac{\left(
-1\right) ^{r}}{\left[ r\right] _{p,q}!}\dfrac{D_{p,q}^{r}\left(
E(-t)\right) }{E(-t)}\text{,}
\end{equation*}
which implies in particular $\left[ p_{n}^{(0)}\right] _{p,q}=\delta _{n}^{0}
$.
\end{definition}

We verify that the calculations of Section 3 transpose to this $p$,$q$%
-analogue, with in particular

\begin{equation*}
D_{p,q}^{r}E(t)=\left[ r\right] _{p,q}!\sum\limits_{n\geq r}\QATOPD[ ] {n}{r}%
_{p,q}e_{n}t^{n-r}\text{,}
\end{equation*}

and

\begin{equation*}
\left[ p_{n}^{(r)}\right] _{p,q}=\left| 
\begin{array}{cccccccc}
\QATOPD[ ] {r}{r}_{p,q}e_{r} & 1 & 0 & 0 & . & . & . & 0 \\ 
\QATOPD[ ] {r+1}{r}_{p,q}e_{r+1} & e_{1} & 1 & 0 & . & . & . & 0 \\ 
\QATOPD[ ] {r+2}{r}_{p,q}e_{r+2} & e_{2} & e_{1} & 1 & 0 & . & . & 0 \\ 
. & . &  & . & . & . & . & . \\ 
. & . &  &  & . & . & . & . \\ 
. & . &  &  &  & . & . & 0 \\ 
. & . &  &  &  &  & . & 1 \\ 
\QATOPD[ ] {n}{r}_{p,q}e_{n} & e_{n-r} & e_{n-r-1} & . & . & . & e_{2} & 
e_{1_{{}}}
\end{array}
\right| \text{.}
\end{equation*}

\textbf{Particular case }$r=1$:\ we also denote $\left[ p_{n}^{(1)}\right]
_{p,q}=\left[ p_{n}\right] _{p,q}$ and we get the following $p,q-$analogues
of the classical case (see [12, Example 8]), with inside the determinants
the $p,q$-analogues of the square brackets:

\begin{equation*}
\left[ p_{n}\right] _{p,q}=\left| 
\begin{array}{ccccccc}
\left[ 1\right] e_{1} & 1 & 0 & 0 & . & . & 0 \\ 
\left[ 2\right] e_{2} & e_{1} & 1 & 0 & . & . & 0 \\ 
. & . & . & . & . & . & . \\ 
. & . & . & . & . & . & . \\ 
. & . & . &  & . & . & 0 \\ 
. & . & . &  &  & . & 1 \\ 
\left[ n\right] e_{n} & e_{n-1} & e_{n-2} & . & . & e_{2} & e_{1}
\end{array}
\right| \text{ ,\ \ \ \ \ \ \ \ \ \ \ \ \ }\left[ n\right]
_{p,q}!\;e_{n}=\left| 
\begin{array}{ccccccc}
\left[ p_{1}\right] & \left[ 1\right] & 0 & 0 & . & . & 0 \\ 
\left[ p_{2}\right] & \left[ p_{1}\right] & \left[ 2\right] & 0 & . & . & 0
\\ 
. & . & . & . & . & . & . \\ 
. & . & . & . & . & . & . \\ 
. & . &  & . & . & . & 0 \\ 
. & . &  &  &  & . & \left[ n-1\right] \\ 
\left[ p_{n}\right] & \left[ p_{n-1}\right] & \left[ p_{n-2}\right] & . & .
& \left[ p_{2}\right] & \left[ p_{1}\right]
\end{array}
\right| \text{.}
\end{equation*}

We leave it to the reader to verify that the proof of $\left( 4.7\right) $\
extends without difficulty to the $p,q$-analog and gives:

\begin{proposition}
We have for all $n\geq r\geq 1$, 
\begin{equation*}
\left[ p_{n}^{(r)}\right] _{p,q}=\sum\limits_{j=r}^{n}p_{n}^{(j)}\sum%
\limits_{l=r}^{j}\left( -1\right) ^{l-r}\binom{j}{l}\QATOPD[ ] {l}{r}_{p,q}%
\text{.}
\end{equation*}
\end{proposition}

Things get more complicated if we want to extend Theorem 4.1. There are
indeed $p,q$-Stirling numbers of the first and second kind [6, 19]. But as
indicated in [6, page 104], there is no $p,q$-analogue of $\left( 4.2\right) 
$ which we used to go from $\left( 4.7\right) $ to $\left( 4.4\right) $.
This apparently prevents getting a $p,q$-analogue of Theorem 4.1.

\section{The polynomials $J_{n}^{\left( r\right) }$ and their linear
recurrence}

We now study the specialization $e_{n}=q^{\binom{n}{2}}/n!$ which
corresponds to the $q$-deformation of the exponential series defined by $%
\left( 1.2\right) $. We begin with two lemmas.

\begin{lemma}
\bigskip \bigskip \bigskip For the $q$-deformation of the exponential series
defined by $\left( 1.2\right) $and $n\geq r+1\geq 2$, we have$:$%
\begin{equation}
p_{n}^{(r)}=\dfrac{\left( 1-q^{r}\right) }{r!}q^{\binom{r}{2}}\left[ p_{n-r}%
\right] _{q^{r}}\text{,}  \tag{6.1}
\end{equation}
where for all $n\geq 1$, $\left[ p_{n}\right] _{q^{r}}$ is deduced from
Equation $\left( 3.6\right) $ with in the determinant $e_{n}=q^{\binom{n}{2}%
}/n!$ and $\left[ n\right] =\left[ n\right] _{q^{r}}$ given by Definition $%
2.1$.
\end{lemma}

\begin{proof}
We start from $\left( 3.4\right) $ with $e_{n}=q^{\binom{n}{2}}/n!$. Then we
factor the first column by $q^{\binom{r}{2}}/r!$ and we subtract the first
column from the second, which gives 
\begin{equation*}
p_{n}^{(r)}=\dfrac{q^{\binom{r}{2}}}{r!}\left| 
\begin{array}{cccccccc}
1 & 0 & 0 & 0 & . & . & . & 0 \\ 
\dfrac{q^{\binom{r+1}{r}-\binom{r}{2}}}{1!} & \dfrac{1-q^{r}}{1!} & 1 & 0 & .
& . & . & 0 \\ 
\dfrac{q^{\binom{r+2}{2}-\binom{r}{2}}}{2!} & q^{\binom{2}{2}}\dfrac{1-q^{2r}%
}{2!} & 1 & 1 & 0 & . & . & 0 \\ 
. & . & . & . & . & . & . & . \\ 
. & . & . & . & . & . & . & . \\ 
\dfrac{q^{\binom{r+j}{2}-\binom{r}{2}}}{j!} & q^{\binom{j}{2}}\dfrac{1-q^{jr}%
}{j!} & . & . & . & . & . & 0 \\ 
. & . & . & . & . & . & . & 1 \\ 
\dfrac{q^{\binom{n}{2}-\binom{r}{2}}}{\left( n-r\right) !} & q^{\binom{n-r}{2%
}}\dfrac{1-q^{\left( n-r\right) r}}{\left( n-r\right) !} & \dfrac{q^{\binom{%
n-r-1}{2}}}{\left( n-r-1\right) !} & . & . & . & \dfrac{q^{\binom{2}{2}}}{2!}
& 1
\end{array}
\right|
\end{equation*}
For $1\leq j\leq n-r$,\ $1-q^{jr}=\left( 1-q^{r}\right) \left[ j\right]
_{q^{r}}$. We therefore factor the second column by $\left( 1-q^{r}\right) $%
, then we expand the determinant with respect to the first row, hence: 
\begin{equation*}
p_{n}^{(r)}=\dfrac{q^{\binom{r}{2}}}{r!}\left( 1-q^{r}\right) \left| 
\begin{array}{ccccccc}
1 & 1 & 0 & . & . & . & 0 \\ 
\left[ 2\right] _{q^{r}}\dfrac{q^{\binom{2}{2}}}{2!} & 1 & 1 & 0 & . & . & 0
\\ 
. & \dfrac{q^{\binom{2}{2}}}{2!} & . & . & . & . & . \\ 
. &  & . & . & . & . & . \\ 
\left[ j\right] _{q^{r}}\dfrac{q^{\binom{j}{2}}}{j!} &  &  & . & . & . & 0
\\ 
. &  &  &  & \;.\; & . & 1 \\ 
\left[ n-r\right] _{q^{r}}\dfrac{q^{\binom{n-r}{2}}}{\left( n-r\right) !} & 
\dfrac{q^{\binom{n-r-1}{2}}}{\left( n-r-1\right) !} & . & . & . & \dfrac{q^{%
\binom{2}{2}}}{2!} & 1
\end{array}
\right|
\end{equation*}
By comparing this determinant to $\left( 3.6\right) $ and taking into
account Definition 2.1 with $\psi (q)=q^{r}$, we obtain the second member of
Equation $\left( 5.1\right) $.
\end{proof}

\bigskip

\begin{lemma}
For the formal series $E_{xq}(t)=\sum\nolimits_{n=0}^{\infty }q^{\binom{n}{2}%
}t^{n}/n!$, we have 
\begin{equation}
D^{r}E_{xq}(t)=q^{\binom{r}{2}}E_{xq}(q^{r}t).  \tag{6.2}
\end{equation}
\end{lemma}

\begin{proof}
On the left handside we have formally\qquad \qquad 
\begin{equation*}
D^{r}E_{xq}(t)=\sum\limits_{n=0}^{\infty }\dfrac{q^{\binom{n}{2}}}{n!}%
D^{r}t^{n}=\sum\limits_{n=r}^{\infty }\dfrac{q^{\binom{n}{2}}}{n!}\dfrac{n!}{%
\left( n-r\right) !}t^{n-r}=\sum\limits_{m=0}^{\infty }q^{\binom{m+r}{2}}%
\dfrac{t^{m}}{m!}\text{.}
\end{equation*}
On the right handside we have\qquad \qquad 
\begin{equation*}
q^{\binom{r}{2}}E_{xq}(q^{r}t)=q^{\binom{r}{2}}\sum\limits_{m=0}^{\infty }q^{%
\binom{m}{2}}\dfrac{\left( q^{r}t\right) ^{m}}{m!}=\sum\limits_{m=0}^{\infty
}q^{\binom{r}{2}+\binom{m}{2}+mr}\dfrac{t^{m}}{m!}\text{,}
\end{equation*}
and we check that 
\begin{equation}
\binom{m+r}{2}=\binom{m}{2}+\binom{r}{2}+mr\text{.}  \tag{6.3}
\end{equation}
\bigskip
\end{proof}

\bigskip It is now possible to introduce very naturally a class of
polynomials with positive integer coefficients.

\begin{theorem}
\bigskip i) For $\ n\geq r\geq 1$ $\ $and $\
E_{xq}(t)=\sum\nolimits_{n=0}^{\infty }q^{\binom{n}{2}}t^{n}/n!$, we have 
\begin{equation}
p_{n}^{(r)}=\left( 1-q\right) ^{n-r}\dfrac{q^{\binom{r}{2}}}{r!\left(
n-r\right) !}J_{n}^{\left( r\right) }(q)\text{,}  \tag{6.4}
\end{equation}
where $J_{n}^{\left( r\right) }$ is a monic polynomial with positive integer
coefficients, a constant term equal to $\left( n-r\right) !$ \ and whose
degree is $\binom{n-1}{2}-\binom{r-1}{2}$. Moreover, \ for all $r\geq 1$ $%
J_{r}^{\left( r\right) }=1$.

\hspace{0.9in}ii) The polynomials $J_{n}^{\left( r\right) }$ satisfy when $%
n-1\geq r\geq 1$, the linear recurrence whose coefficients belong to $%
\mathbb{N}\left[ q\right] $: 
\begin{equation}
J_{n}^{\left( r\right) }(q)=\sum\limits_{j=1}^{n-r}\left[ r\right]
_{q}^{j}\;q^{\binom{j}{2}}\binom{n-r}{j}J_{n-r}^{\left( j\right) }(q)\text{.}
\tag{6.5}
\end{equation}
This recurrence suffices with the initial conditions $J_{r}^{\left( r\right)
}=1$ for $r\geq 1$, to calculate all the $J_{n}^{\left( r\right) }$ for $%
n\geq r+1$.
\end{theorem}

\begin{proof}
\textbf{Let first prove}\textit{\ i)}

a) For $n=r.$ From $\left( 3.2\right) $ and Lemma 6.2 we have 
\begin{equation*}
\sum\limits_{n\geq r}p_{n}^{(r)}\left( -t\right) ^{n-r}=\dfrac{1}{r!}\dfrac{%
D^{r}E_{xq}(t)}{E_{xq}(t)}=\dfrac{q^{\binom{r}{2}}}{r!}\dfrac{E_{xq}(q^{r}t)%
}{E_{xq}(t)}\text{.}
\end{equation*}

The constant term of the formal series quotient $E_{xq}(q^{r}t)/E_{xq}(t)$
is $E_{xq}(q^{r}0)/E_{xq}(0)=1$. Thefore the constant term of the series on
the left handside of the above equation is $p_{r}^{(r)}=q^{\binom{r}{2}}/r!$%
. By comparison to $\left( 6.4\right) $ we get \ $J_{r,r}=1$ \ and this
polynomial verifies all the assertions of item \textit{\ i)} of the theorem.

b) For the general case $n\geq r$ $\geq 1$ we reason by induction on $n$.
For $n=r=1$ it is a particular case of a) and it is therefore proved. Let's
assume true the theorem for all pairs of integers $\left( l,j\right) $ such
that $1\leq l$ $\leq n$ and $1\leq j\leq l$. And let's prove it for $n+1$
and all $r$ between $1$ and $n$ (since it's already proved for $r=n+1$ by
a)). We have with Lemma 6.1 
\begin{equation}
p_{n+1}^{(r)}=\left( 1-q^{r}\right) \dfrac{q^{\binom{r}{2}}}{r!}\left[
p_{n+1-r}\right] _{q^{r}}\text{,}  \tag{6.6}
\end{equation}
and the particular case $\left( 4.8\right) $ of Theorem 4.1 gives 
\begin{equation*}
\left[ p_{n+1-r}\right] _{q^{r}}=\sum\limits_{j=1}^{n+1-r}\left(
1-q^{r}\right) ^{j-1}p_{n+1-r}^{(j)}\text{.}
\end{equation*}
But for $1\leq r\leq n$ we have $1\leq n+1-r\leq n$, so the induction
hypothesis applies to $p_{n+1-r}^{(j)}$, hence

\begin{equation*}
\left[ p_{n+1-r}\right] _{q^{r}}=\sum\limits_{j=1}^{n+1-r}\left(
1-q^{r}\right) ^{j-1}\left( 1-q\right) ^{n+1-r-j}\dfrac{q^{\binom{j}{2}}}{%
j!\left( n+1-r-j\right) !}J_{n+1-r}^{\left( j\right) }\text{,}
\end{equation*}
by replacing this expression in ($6.6)$ we get after some easy
transformations \ 
\begin{equation*}
p_{n+1}^{(r)}=\left( 1-q\right) ^{n+1-r}\dfrac{q^{\binom{r}{2}}}{r!\left(
n+1-r\right) !}\sum\limits_{j=1}^{n+1-r}\left[ r\right] ^{j}\binom{n+1-r}{j}%
q^{\binom{j}{2}}J_{n+1-r}^{\left( j\right) }\text{.}
\end{equation*}

Let 
\begin{equation}
J_{n+1}^{\left( r\right) }(q)=\sum\limits_{j=1}^{n+1-r}\left[ r\right] ^{j}%
\binom{n+1-r}{j}q^{\binom{j}{2}}J_{n+1-r}^{\left( j\right) }\text{,} 
\tag{6.7}
\end{equation}
then we have verified $\left( 6.4\right) $ for $n+1$. Moreover, since $%
\left( 6.7\right) $ links $J_{n+1}^{\left( r\right) }$ to polynomials which
have positive integer coefficients by the induction hypothesis with
coefficients which are themselves polynomials in q with positive integer, it
is the same for $J_{n+1}^{\left( r\right) }$.

On the right hand side of $\left( 6.7\right) $ the degree of each term of
the sum is by the induction hypothesis 
\begin{equation*}
d%
{{}^\circ}%
\left( \left[ r\right] ^{j}\binom{n+1-r}{j}q^{\binom{j}{2}}J_{n+1-r}^{\left(
j\right) }\right) =(r-1)j+\binom{j}{2}+\binom{n-r}{2}-\binom{j-1}{2}=\binom{%
n-r}{2}+rj-1\text{.}
\end{equation*}
It is clearly maximum for $j=n+1-r$ and is then $\binom{n}{2}-\binom{r-1}{2}$
(use $\left( 6.3\right) $ with $m=n-r$). The coefficient of the monomial of
maximum degree $n+1-r$, is $J_{n+1-r}^{\left( n+1-r\right) }$ which is equal
to 1 by a). $J_{n+1}^{\left( r\right) }$ is therefore a monic \ polynomial.

When $j=1$ the term of the sum is $\left[ r\right] \left( n+1-r\right)
J_{n+1-r}^{\left( 1\right) }(q)$. This polynomial has an order equal to zero
and it is the only term of the sum having this property. The order of $%
J_{n+1}^{\left( r\right) }$ is therefore zero and its constant term is by
the induction hypothesis $\left( n+1-r\right) .\left( n-r\right) !=\left(
n+1-r\right) !.$ Which ends the proof of \textit{i)}.

\textbf{Proof of item} \textit{ii)} Equation $\left( 6.5\right) $ is none
other than the equation $\left( 6.7\right) $ defining $J_{n+1}^{\left(
r\right) }$ in the proof of \textit{i)}. $\left( 6.5\right) $ is therefore
proved for all $n$ and $r$ such that $n-1\geq r\geq 1$. Arrange the
polynomials in a triangular table as shown below (we can also set $%
J_{n}^{\left( r\right) }=0$ if $n<r$ by convention). We see that linear
horizontal recurrence $\left( 6.5\right) $ makes it possible to determine
row by row all the polynomials such as $n>r$.
\end{proof}

\begin{tabular}{|l|l|l|l|l|l|}
\hline
\begin{tabular}{ll}
& $r$ \\ 
$n$ & 
\end{tabular}
& $1$ & $2$ & $3$ & $4$ & $5$ \\ \hline
$1$ & $1$ &  &  &  &  \\ \hline
$2$ & $1$ & $1$ &  &  &  \\ \hline
$3$ & $2+q$ & $1+q$ & $1$ &  &  \\ \hline
$4$ & $6+6q+3q^{2}+q^{3}$ & $2+3q+2q^{2}+q^{3}$ & $1+q+q^{2}$ & $1$ &  \\ 
\hline
$5$ & $
\begin{tabular}{l}
$24+36q+30q^{2}+20q^{3}$ \\ 
$+10q^{4}+4q^{5}+q^{6}$%
\end{tabular}
$ & $
\begin{tabular}{l}
$6+12q+12q^{2}+10q^{3}$ \\ 
$+6q^{4}+3q^{5}+q^{6}$%
\end{tabular}
$ & 
\begin{tabular}{l}
$2+3q+4q^{2}+$ \\ 
$3q^{3}+2q^{4}+q^{5}$%
\end{tabular}
& $1+q+q^{2}+q^{3}$ & $1$ \\ \hline
\end{tabular}

\medskip \textbf{Table 1: Polynomials }$J_{n}^{\left( r\right) }$\textbf{\
for }$5\geq n\geq r\geq 1.$

\medskip \textbf{Exemple of application of Recurrence }$\left( 6.5\right) $%
\textbf{\ for} $n=6,r=2$:$\qquad $%
\begin{equation*}
J_{6}^{\left( 2\right) }=4\left( 1+q\right) J_{4}^{\left( 1\right) }+6\left(
1+q\right) ^{2}qJ_{4}^{\left( 2\right) }+4\left( 1+q\right) ^{3}q^{3}J_{4}^{
\left[ 3\right] }+\left( 1+q\right) ^{4}q^{6}\text{,}
\end{equation*}

hence:$\qquad \qquad J_{6}^{\left( 2\right)
}=24+60q+78q^{2}+80q^{3}+68q^{4}+52q^{5}+35q^{6}+20q^{7}+10q^{8}+4q^{9}+q^{10} 
$.

\bigskip

\textbf{Case} $r=1$: $\left( 6.5\right) $ gives for all $n\geq 1$, $\qquad
J_{n+1}^{\left( 1\right) }(q)=\sum\limits_{j=1}^{n}\binom{n}{j}q^{\binom{j}{2%
}}J_{n}^{\left( j\right) }(q)$.

Recurrence $\left( 6.5\right) $\ can be generalized in the following form,
by setting for all $n\geq 0$, $J_{n}^{\left( 0\right) }=\delta _{n,0}$ and $%
\left[ 0\right] _{q}^{0}=1$.

\begin{corollary}
\textit{We have for all pair of integers }$\left( n,r\right) $\textit{\ such
that }$n\geq r\geq 0$\textit{\ } 
\begin{equation}
J_{n}^{\left( r\right) }(q)=\sum\limits_{j=0}^{n-r}\left[ r\right]
_{q}^{j}\;q^{\binom{j}{2}}\binom{n-r}{j}J_{n-r}^{\left( j\right) }(q)\text{.}
\tag{6.8}
\end{equation}
\end{corollary}

\begin{proof}
If \ $n>$ $r\geq 1$,$\;\;\left( 6.5\right) $ implies $\left( 6.8\right) $
since the added term contains $J_{n-r}^{\left( 0\right) }$ which is zero.

If $n=r\geq 1$, the left side of $\left( 6.8\right) $ is $1$ and the right
side is $\left[ 0\right] ^{0}q^{0}\binom{0}{0}J_{0}^{\left( 0\right) }=1.$

If $n\geq r=0$, the left side of $\left( 6.8\right) $ is $J_{n}^{\left(
0\right) }=\delta _{n,0}$ and on the right side, the sum is reduced to the
first term $\left[ 0\right] ^{0}q^{0}\binom{0}{0}J_{n}^{\left( 0\right)
}=\delta _{n,0}$.
\end{proof}

Now let's make the connection with the polynomials introduced previously of
which we briefly remind the genesis. We refer to [21] for more details. In
[13], Mallow and Riordan defined the inversion polynomial enumerator for
rooted trees with $n$ nodes, that they denoted $J_{n}$. Then Stanley $\left[
16\right] $ and Yan $\left[ 20\right] $ successively generalized the
definition of inversion polynomial enumerator to sequences of \ ''colored''
rooted forests, which correspond to classical parking functions in the
terminology of [21, Section 1.4.4]. These polynomials are characterized by a
pair of parameters $\left( a,b\right) $ and denoted $I_{m}^{(a,b)}$. It has
been proven combinatorially that [20, Corollary 5.1]: 
\begin{equation}
\sum\limits_{m\geq 0}\left( q-1\right) ^{m}I_{m}^{\left( a,b\right) }\dfrac{%
t^{m}}{m!}=\dfrac{\sum\limits_{m\geq 0}q^{am+b\binom{m}{2}}\dfrac{t^{m}}{m!}%
}{\sum\limits_{m\geq 0}q^{b\binom{m}{2}}\dfrac{t^{m}}{m!}}\text{.}  \tag{6.9}
\end{equation}

Note that these polynomials have other combinatorial representations and are
still the subject of recent research (see for example [10, 18]).

\begin{corollary}
We have for all $n\geq r\geq 1$, 
\begin{equation}
J_{n}^{\left( r\right) }=I_{n-r}^{\left( r,1\right) }\Leftrightarrow
I_{m}^{\left( r,1\right) }=J_{m+r}^{\left( r\right) }\text{.}  \tag{6.10}
\end{equation}
\end{corollary}

\begin{proof}
\bigskip For $a=r,b=1$, $\left( 6.9\right) $ becomes 
\begin{equation}
\sum\limits_{m\geq 0}\left( q-1\right) ^{m}I_{m}^{\left( r,1\right) }\dfrac{%
t^{m}}{\left( m\right) !}=\dfrac{\sum\limits_{m\geq 0}q^{rm+\binom{m}{2}}%
\dfrac{t^{m}}{\left( m\right) !}}{\sum\limits_{m\geq 0}q^{\binom{m}{2}}%
\dfrac{t^{n}}{n!}}\text{.}  \tag{6.11}
\end{equation}
The denominator of this fraction is $E_{xq}\left( t\right) $. For the
numerator we have by Lemma 2, $D^{r}E_{xq}\left( t\right) =q^{\binom{r}{2}%
}E_{xq}\left( q^{r}t\right) =q^{\binom{r}{2}}\sum\nolimits_{m\geq 0}q^{rm+%
\binom{m}{2}}t^{m}/m!$. So by multiplying $\left( 6.11\right) $ by $q^{%
\binom{r}{2}}/r!$ and with $m=n-r$, 
\begin{equation}
\sum\limits_{n\geq r}\left( q-1\right) ^{n-r}q^{\binom{r}{2}}I_{n-r}^{\left(
r,1\right) }\dfrac{t^{n-r}}{r!\left( n-r\right) !}=\dfrac{1}{r!}\dfrac{%
D^{r}E_{xq}\left( t\right) }{E_{xq}\left( t\right) }\text{.}  \tag{6.12}
\end{equation}
By identifying the coefficients of $\left( 6.12\right) $ and those of $%
\left( 3.2\right) $ we therefore obtain for $E_{xq}\left( t\right) $%
\begin{equation*}
p_{n}^{\left( r\right) }=\left( 1-q\right) ^{n-r}\dfrac{q^{\binom{r}{2}}}{%
r!\left( n-r\right) !}I_{n-r}^{\left( r,1\right) }\text{,}
\end{equation*}
that is to say by comparison with $\left( 6.4\right) $,\ $J_{n}^{\left(
r\right) }=I_{n-r}^{\left( r,1\right) }$.
\end{proof}

\bigskip Case $r=1$ corresponds to the case of a tree [13]. So we have $%
J_{n}^{\left( 1\right) }=J_{n}$, which gives another notation for $%
J_{n}^{\left( 1\right) }$.

We observe that our study clearly differs from those in the previously cited
publications. In those, the polynomials $I_{m}^{\left( a,b\right) }$\ are
introduced in a combinatorial manner, and Equation $\left( 6.9\right) $\ is
also derived combinatorially. In contrast, we have formally define the
polynomials $J_{n}^{\left( r\right) }$ through a recursive application of
Theorem 4.1. Recurrence $\left( 6.5\right) $, inherent in the proof of
Theorem 6.3, will in turn lead to a new combinatorial representation of the $%
J_{n}^{\left( r\right) }$ in Section 9. Note, however, that our approach
applies only to the class of polynomials given by $\left( 6.10\right) $.

\section{\protect\bigskip Reciprocal polynomials}

\begin{definition}
\bigskip The polynomials $\overline{J_{n}^{\left( r\right) }}$ are defined
by 
\begin{equation}
\overline{J_{n}^{\left( r\right) }}\left( q\right) =q^{\binom{n-1}{2}-\binom{%
r-1}{2}}J_{n}^{\left( r\right) }(1/q)\text{.}  \tag{7.1}
\end{equation}
These are the reciprocal polynomials of $J_{n}^{\left( r\right) }$
\end{definition}

Since $\overline{\overline{J_{n}^{\left( r\right) }}}=J_{n}^{\left( r\right)
}$, $\left( 7.1\right) $ is equivalent to 
\begin{equation}
J_{n}^{\left( r\right) }\left( q\right) =q^{\binom{n-1}{2}-\binom{r-1}{2}}%
\overline{J_{n}^{\left( r\right) }}(1/q)\text{.}  \tag{7.1bis}
\end{equation}
It is well known that these reciprocal polynomials are equal to the sum
enumerator of the corresponding parking functions which definition we
recall. We will refer for more details to [21] from which we almost adopt
the notations for the parking functions of our case, i.e when parameters $%
\left( a,b\right) =\left( r,1\right) $. Let $m=n-r$ and let $PK_{m}(r,1)$
designate the set of parking functions associated with the sequence $%
r,r+1,r+2,...,r+m-1$. $PK_{m}(r,1)$ is the set of sequences of integers $%
\mathbf{a}=\left( a_{1},a_{2},...,a_{m}\right) $, such that : 
\begin{equation*}
a_{(i)}<r+i-1\text{ \ \ for \ \ }1\leq i\leq m\text{,}
\end{equation*}
$\left( a_{(1)},a_{\left( 2\right) },...,a_{\left( m\right) }\right) $ being
the sequences of $a_{i}$ ordered in a non-decreasing way.

We therefore have with $\left| \mathbf{a}\right| =a_{1}+a_{2}+...+a_{m}$ 
\begin{equation*}
\overline{J_{m+r}^{\left( r\right) }}(q)=\sum\limits_{\mathbf{a}\in
PK_{m}(r,1)}q^{\left| \mathbf{a}\right| }=S_{m}\left( q;r\right) \text{,}
\end{equation*}
where the second member is the sum enumerator written $S_{m}\left(
q;r\right) $ (instead of $S_{m}\left( q;r,r+1,r+2,...,r+m-1\right) $ in
[21]). We verify that this relation is still true if $m=0$, on the condition
that we set $PK_{0}(r,1)=\emptyset $ and $\sum\limits_{\mathbf{a}\in
\emptyset }q^{\left| \mathbf{a}\right| }=1$.

With $\left( 7.1bis\right) $, $\left( 6.5\right) $ gives:

\begin{corollary}
The reciprocal polynomials $\overline{J_{n}^{(r)}}$ \ satisfy when $n-1\geq
r\geq 1$, the linear recurrence of which coefficients are elements of $%
\mathbb{N}[q]$ 
\begin{equation}
\text{ \ \ \ }\overline{J_{n}^{\left( r\right) }}(q)=\sum\limits_{j=1}^{n-r}%
\left[ r\right] _{q}^{j}\;q^{r\left( n-r-j\right) }\binom{n-r}{j}\overline{%
J_{n-r}^{\left( j\right) }}(q)  \tag{7.2}
\end{equation}
\end{corollary}

We could alternatively use this recurrence to find the table of $\overline{%
J_{n}^{\left( r\right) }}$ with $\overline{J_{r}^{\left( r\right) }}%
=1.\bigskip $

In [11] Kung and Yan presented an alternative study of parking functions
from Goncarov polynomials. It is interesting to compare our linear
recurrence $\left( 7.2\right) $ with a linear recurrence resulting from this
study. This is Equation $\left( 6.2\right) $\ in [11] (or $\left(
1.27\right) $ in [21]), which in our case and with our notations is written: 
\begin{equation*}
1=\sum\limits_{k=0}^{m}\binom{m}{k}q^{\left( r+k\right) \left( m-k\right)
}\left( 1-q\right) ^{k}\overline{J_{k+r}^{\left( r\right) }}\text{,}
\end{equation*}
which is equivalent with $m+r=n$ \ and $l=k+r$,\ to 
\begin{equation}
\left( 1-q\right) ^{n-r}\overline{J_{n}^{\left( r\right) }}%
(q)=1-\sum\limits_{l=r}^{n-1}\binom{n-r}{l-r}q^{l\left( n-l\right) }\left(
1-q\right) ^{l-r}\overline{J_{l}^{\left( r\right) }}\left( q\right) \text{.}
\tag{7.3}
\end{equation}
To see the difference, let's take the example of $n=6,r=2$:

\begin{equation}
\overline{J_{6}^{\left( 2\right) }}=4\left( 1+q\right) q^{6}\overline{%
J_{4}^{\left( 1\right) }}+6\left( 1+q\right) ^{2}q^{4}\overline{%
J_{4}^{\left( 2\right) }}+4\left( 1+q\right) ^{3}q^{2}\overline{%
J_{4}^{\left( 3\right) }}+\left( 1+q\right) ^{4}\overline{J_{4}^{\left(
4\right) }}\text{,}  \tag{with 7.2}
\end{equation}

\begin{equation}
\left( 1-q\right) ^{4}\overline{J_{6}^{\left( 2\right) }}=1-q^{8}\overline{%
J_{2}^{\left( 2\right) }}-4q^{9}\left( 1-q\right) \overline{J_{3}^{\left(
2\right) }}-6q^{8}\left( 1-q\right) ^{2}\overline{J_{4}^{\left( 2\right) }}%
-4q^{5}\left( 1-q\right) ^{3}\overline{J_{5}^{\left( 2\right) }}\text{,} 
\tag{with 7.3}
\end{equation}

Let call Table 2, the table - not shown here - deduced from Table 1 by
replacing the polynomials $J_{n}^{\left( r\right) }$ by their reciprocals $%
\overline{J_{n}^{\left( r\right) }}$. The two equations above show that:

* On the one hand, our recurrence relation $\left( 7.2\right) $ is
horizontal using here row 4 of Table\ $2$, and the coefficients belong to $%
\mathbb{N}\left[ q\right] $.

* On the other hand, the recurrence relation $\left( 7.3\right) $ is
vertical using here column 2 of Table $2$, and the coefficients belong to $%
\mathbb{Z}\left[ q\right] $.

\section{\protect\bigskip An explicit formula for the $J_{n}^{(r)}$}

From the recurrence formulas of Section 6, we will prove the following
theorem which is declined in two versions corresponding to the two versions $%
\left( 6.5\right) $ and $\left( 6.8\right) $ of the recurrence relation. We
first generalize the notations used for integer partitions.

\begin{notation}
Let $u=\left( u_{1},u_{2},...,u_{k}\right) $ be a sequence of strictly
positive integers, we set 
\begin{equation*}
\left| u\right| =u_{1}+u_{2}+...+u_{k}\text{ \ \ \ and \ \ }n(u^{\prime })=%
\binom{u_{1}}{2}+\binom{u_{2}}{2}+...+\binom{u_{k}}{2}\text{.}
\end{equation*}
\end{notation}

Note that $u^{\prime }$\ is not defined in this case, only $n\left(
u^{\prime }\right) $ has a meaning given by the above formula. More
generally if $\mathcal{U}$\ is the set of infinite sequences of integers $%
u=\left( u_{1},u_{2},...\right) $ with $u_{i}\geq 0$, such that only a
finite number of the term $u_{i}$ are non-zero, we can extend the above
notations to all $u\in \mathcal{U}$ by

\begin{equation*}
\left| u\right| =\sum\limits_{i\geq 1}u_{i}\text{ \ \ and \ }n(u^{\prime
})=\sum\limits_{i\geq 1}\binom{u_{i}}{2}\text{.}
\end{equation*}

\begin{theorem}
a) For any pair of integers $\left( n,r\right) $ satisfying $n-1\geq r\geq 1$
\ we have the explicit formula 
\begin{equation}
J_{n}^{\left( r\right) }(q)=\sum \left[ r\right] _{q}^{u_{1}}\left[ u_{1}%
\right] _{q}^{u_{2}}...\left[ u_{k-1}\right] _{q}^{u_{k}}q^{n(u^{\prime })}%
\binom{n-r}{u_{1},u_{2},...,u_{k}}\text{,}  \tag{8.1}
\end{equation}
where the sum is over the $k$-multiplets of strictly positive integers $%
u=\left( u_{1},u_{2},...,u_{k}\right) $ with $k\geq 1$ and $\left| u\right|
=n-r$, with $n(u^{\prime })=\sum\limits_{i=1}^{k}\binom{u_{i}}{2}$\ and
where $\binom{n-r}{u_{1},u_{2},...,u_{k}}$ is the multinomial coefficient.

\hspace{0.8in}b) For any pair of integers $\left( n,r\right) $ satisfying $%
n\geq r\geq 0$ \ we have: 
\begin{equation}
J_{n}^{\left( r\right) }(q)=\left( n-r\right) !\sum \left[ r\right]
_{q}^{u_{1}}q^{n(u^{\prime })}\prod\limits_{i\geq 1}\dfrac{\left[ u_{i}%
\right] _{q}^{u_{i+1}}}{u_{i}!}\text{,}  \tag{8.2}
\end{equation}
where the sum is over the sequences $u=\left( u_{i}\right) _{i\geq 1}$ $\in $
$\mathcal{U}$ satisfying $\left| u\right| =\sum\nolimits_{i\geq 1}u_{i}=n-r$
and with $n(u^{\prime })=\sum\nolimits_{i\geq 1}\binom{u_{i}}{2}$.
\end{theorem}

As we will see below, the sum in $\left( 8.2\right) $\ contains only a
finite number of non-zero terms. Note that if b) is apparently more
complicated than a), it is in fact easier to prove.\medskip

\begin{proof}
We first show the equivalence of a) and b), when $n-1\geq r\geq 1.$

In order for the general term of the sum in (8.2) to be non-zero, it is
necessary - taking into account $\left[ 0\right] ^{n}=\delta _{n}^{0}$ -
that the finite number (say $k$) of non-zero terms of the sequence $u$,
occupy the first $k$ places in this sequence. In other words, the summation
can be limited to these sequences, which we call commencing sequences, and
we denote by $\left( 8.2bis\right) $ the summation deduced from $\left(
8.2\right) $\ and coresponding to these commencing sequences. The
application, which associates to a multiplet of the sum of (8.1) $\left(
u_{1},u_{2},...,u_{k}\right) $ with $u_{1}+u_{2}+...+u_{k}=n-r$ $\ $and $%
k\geq 1$, the commencing sequence $u=\left(
u_{1},u_{2},...,u_{k},0,0,...\right) $ with $\left| u\right| =n-r$ ,\ is
clearly a bijection between the two summation sets of $\left( 8.1\right) $
and $\left( 8.2bis\right) $. Moreover, we check that the value of the terms
in the respective sums, thus put in 1-1 corespondence are equal. The three
sums $\left( 8.1\right) $, $\left( 8.2\right) $ and $\left( 8.2bis\right) $\
are therefore equal.

Let now prove $\left( 8.2bis\right) $ for all $n\geq r\geq 0$.

1) If $n=r$, the condition $\left| u\right| =\sum\nolimits_{i\geq
1}u_{i}=n-r=0$ implies $u_{i}=0$ for all $i\geq 1$. The sum in $\left(
8.2bis\right) $ therefore reduces to $0!\left[ r\right] ^{0}q^{0}\left[ 0%
\right] ^{0}.../0!..=1$ which is indeed equal to $J_{r}^{\left( r\right) }$
for all $r\geq 0$.

2) If $r=0$ and $n>0$, it is necessary for a commencing sequence $u$
verifying $\left| u\right| =n-r>0$, that its first term $u_{1}$ be non-zero.
Therefore, the term indexed by $u$ in the sum of ($8.2bis)$ is zero, since
its contains the factor $\left[ 0\right] ^{u_{1}}=0$. Hence, the sum in ($%
8.2bis$) is zero, which is indeed $J_{n}^{\left( 0\right) }$ for $n>0$.

3) Let now prove by induction on $n$ the general case $n\geq r\geq 0$.

When $n=r=0$, it has already been checked by 1). Assume that relation ($%
8.2bis$) is true for any pair of integers $\left( l,j\right) $ such that $%
0\leq l\leq n-1$, \ $0\leq j\leq l$ and show that it is true for $l=n$, $n>0$
and $0\leq j=r\leq n$.

If $r=0$, it is the case 2) already proved since $n>0$.

If $r=n$, it is the case 1) already proved.

So suppose $1\leq r\leq n-1$. The set of commencing sequences verifying $%
\left| u\right| =\sum\nolimits_{i\geq 0}u_{i}=n-r>0$ can be decomposed
according to the value of $u_{1}$ into:

* $u_{1}=1$, from where $\sum\nolimits_{i\geq 2}u_{i}=n-r-1$,

and in a general way for $1\leq j\leq n-r$:

* $u_{1}=j$,\ from where $\sum\nolimits_{i\geq 2}u_{i}=n-r-j$.

Equation $\left( 8.2bis\right) $\ can therefore be written 
\begin{equation}
J_{n}^{\left( r\right) }\left( q\right) =\sum\limits_{j=1}^{n-r}\left[ r%
\right] ^{j}q^{\binom{j}{2}}\binom{n-r}{j}\left[ \left( n-r-j\right)
!\sum\limits_{u_{2}+u_{3}+...=n-r-j}\left[ j\right] ^{u_{2}}q^{\binom{u_{2}}{%
2}+\binom{u_{3}}{2}+...}\prod\limits_{i\geq 2}\dfrac{\left[ u_{i}\right]
^{u_{i+1}}}{u_{i}!}\right] \text{.}  \tag{8.3}
\end{equation}
Let $v_{i}=u_{i+1}$ for $i\geq 1$, then\ \ $v=\left( v_{i}\right) _{i\geq 1}$
is a sequence satisfying $\left| v\right| =(n-r)-j$ and the square bracket
in $\left( 8.3\right) $ can be rewritten 
\begin{equation*}
\left( (n-r)-j\right) !\sum \left[ j\right] _{q}^{v_{1}}q^{n(v^{\prime
})}\prod\limits_{i\geq 1}\dfrac{\left[ v_{i}\right] _{q}^{v_{i+1}}}{v_{i}!}%
\text{,}
\end{equation*}
where the sum is over sequences $\nu $ of integers such that $\left|
v\right| =n-r-j$. But this sum is $J_{n-r}^{\left( j\right) }$ according to
the induction hypothesis.The second member of $\left( 8.3\right) $ is
therefore 
\begin{equation*}
\sum\limits_{j=1}^{n-r}\left[ r\right] _{q}^{j}\;q^{\binom{j}{2}}\binom{n-r}{%
j}J_{n-r}^{\left( j\right) }(q\text{,}
\end{equation*}
which is indeed equal to $J_{n}^{\left( r\right) }$, according to $\left(
6.5\right) $.
\end{proof}

.

\bigskip \textbf{Case }$q=1$. We obtain with $\left( 8.1\right) $ 
\begin{equation}
J_{n}^{\left( r\right) }(1)=\sum_{\substack{ u_{1}+u_{2}+...+u_{k}=n-r  \\ %
u_{i}\geq 1,\;k\geq 1}}r^{u_{1}}u_{1}^{u_{2}}...u_{k-1}^{u_{k}}\binom{n-r}{%
u_{1},u_{2},...,u_{k}}=rn^{n-r-1}\text{.}  \tag{8.4}
\end{equation}

It follows from the generalization, made in [20] to interpret polynomials $%
I_{m}^{(a,b)}$ as inversion enumerators, that $J_{n}^{\left( r\right)
}\left( 1\right) $ is equal to the number of forests on $n$ vertices
(including the roots), comprising $r$ rooted trees whith specified roots.
The set of these forests is denoted by $\mathcal{F}_{m}\left( r,1\right) $
in [21] with here $m=n-r$. Note that according to [17, Proposition 5.3.2],
the number of these forests is $J_{n}^{\left( r\right) }(1)=rn^{n-r-1}$,
which can also easily be verified by induction.

Equation $\left( 8.4\right) $\ can be linked to a formula given by Katz,
whose description and reference can be found in [14, page 19]. This formula
gives the number of connected directed graph or functionnal digraph, made up
of rooted trees whose roots determine a directed cycle. If the cycle has
length $r$ and if the number of nodes including the roots is $n$, then this
number is with our notations [14, Formula page 19]:

\begin{equation}
D\left( n,r\right) =\left( r-1\right) !\sum_{\substack{ %
u_{1}+u_{2}+...+u_{k}=n-r  \\ u_{i}\geq 1,\;k\geq 1}}%
r^{u_{1}}u_{1}^{u_{2}}...u_{k-1}^{u_{k}}\binom{n-r}{u_{1},u_{2},...,u_{k}}%
\text{.}  \tag{8.5}
\end{equation}

It is easy to deduce $\left( 8.5\right) $ from $\left( 8.4\right) $ by
multiplying the number $J_{n}^{\left( r\right) }(1)$ of forests with $r$
rooted trees by the number of cycles $\left( r-1\right) !$, that can be
formed with the roots. On the other hand, Katz's proof of $\left( 8.5\right) 
$ consists to enumerate the digraphs according to the numbers $u_{i}$ of
nodes located at distance $i$ from the cycle. This suggests that conversely,
Equation $\left( 8.1\right) $ could be interpreted, putting aside the factor 
$\left( r-1\right) !$, as a $q$-refinement of the Katz's enumeration. It is
this statistic that is the subject of Section 9.

\section{\protect\bigskip Level statistics on forests}

We denote $\mathcal{F}_{n,R}$ the set of rooted forests whose vertices are $%
V=\mathbf{n}=\left\{ 1,2,...,n\right\} $, and whose roots are a specified
subset $R\subseteq \mathbf{n}$ with $\left| R\right| =r$. $\mathcal{F}_{n,R}$
is equipotent to $\mathcal{F}_{n-r}\left( r,1\right) $ of [21]. We assume
until further notice that $n-1\geq r\geq 1$, therefore $R\subset \mathbf{n}$%
. Let $F\in \mathcal{F}_{n,R}$, $k$ will designate the height of $\ F$, i.e.
the greatest height of its trees. For $i$ an integer between $0$ and $k$,
let $V_{i}=\left\{ v\in \mathbf{n}\text{;}\;\text{distance of }v\text{ from
the root}=i\right\} $ and $u_{i}=\left| V_{i}\right| $; in particular $%
V_{0}=R$ and $u_{0}=r$. The map which to $v\in \mathbf{n}$ associates $V_{i}$
such that $v\in V_{i}$\ \ is denoted $L_{F}$, $L_{F}\left( v\right) $ is the
level (or generation) of $v$.

We necessarily have $u_{0}+u_{1}+....u_{k}=n$ . According to the previous
section, we have with these notations and $r\leq n-1$ which implies $k\geq 1$

\begin{equation*}
\left| \mathcal{F}_{n,R}\right| =J_{n}^{\left( r\right) }\left( 1\right)
=\sum_{\substack{ u_{0}+u_{1}+...+u_{k}=n  \\ u_{0}=r,\;u_{i}\;\geq 1}}%
u_{0}^{u_{1}}u_{1}^{u_{2}}...u_{k}^{u_{k+1}}\binom{n-r}{u_{1},u_{2},...,u_{k}%
}\text{.}
\end{equation*}

Level statistics can be described quite generally with the following
definitions.

\begin{definition}
Let\textbf{\ }$\mathbf{2}^{\mathbf{n}}$ be the set of subsets of $\mathbf{n}$%
. A ranking $\rho $\ on $\mathbf{2}^{n}$ \ is the data for each $P\in 
\mathbf{2}^{n}$ $\ $of a bijection from $P$ into $\mathbf{p}=\left\{
1,...,p\right\} $ where $p=\left| P\right| $, bijection denoted $\rho \left(
P\right) $.
\end{definition}

\textbf{Exemples of ranking on }$\mathbf{2}^{\mathbf{n}}$: the increasing
ranking $\rho _{+}$ is the one for which $\rho _{+}\left( P\right) $ is the
increasing bijection for all $P\in \mathbf{2}^{\mathbf{n}}$; the decreasing
ranking $\rho _{-}$ is the one for which $\rho _{-}\left( P\right) $ is the
decreasing bijection for all $P\in \mathbf{2}^{\mathbf{n}}$. Any combination
of $\rho _{+}$ and $\rho _{-}$ is also a ranking on\textbf{\ }$\mathbf{2}^{%
\mathbf{n}}$, for example that which is equal to $\rho _{+}$ if $\left|
P\right| $ is even and to $\rho _{-}$ otherwise.

\begin{definition}
\bigskip\ Let $\rho $ be a ranking on\textbf{\ }$\mathbf{2}^{\mathbf{n}}$
and $\mathcal{F}_{\mathbf{n},R}$ the set of forests defined above, $n$ and $%
R $ being given. The weight associated with $\rho $ of a vertex $v$ of $F\in 
\mathcal{F}_{n,R}$ is $w_{\rho }(v)=\rho \left( L_{F}\left( v\right) \right)
\left( v\right) $
\end{definition}

Note that this weight depends only on $\rho $ and the generation of $v$. We
note $p\left( v\right) $ the (only) parent of the vertex $v\in \mathbf{n}-R$.

\bigskip

\textbf{Example of a weighted forest with }$n=13$\textbf{, }$r=3,$ $k=4$, $%
R=\left\{ 7,11,2\right\} $\textbf{\ and }$\rho _{+}$. The forest $F\in 
\mathcal{F}_{13,R}$ is represented in Figure 1. For each of its vertices $v$%
\ labeled in black, we have indicated in red the weight associated with $%
\rho _{+}$. The cardinal of each of its level is also indicated on the right.

\begin{figure}[tbph]
\centering\includegraphics[width=14cm]{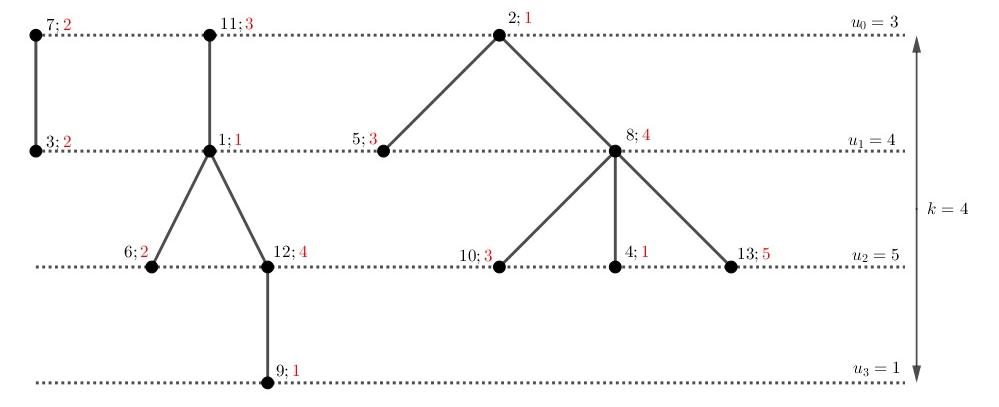}
\caption[Exemple of a weighted forest ]{A weighted forest $F\in \mathcal{F}%
_{13,R}$ with $R=\left\{ 7,11,2\right\} $. The label of each vertex is in
black, its weight is in red. $u_{i}=\left| V_{i}\right| $ for $0\leq i\leq
3. $}
\end{figure}

\bigskip

\begin{definition}
Let $\rho $ a ranking on $\mathbf{2}^{n}$, for all $F\in \mathcal{F}_{n,R}$ (%
$n$ and $R\subset \mathbf{n}$ given), the level statistic $l_{\rho }$
associated with $\rho $ is defined by 
\begin{equation*}
l_{\rho }\left( F\right) =n(u^{\prime })+\sum\nolimits_{v\in \mathbf{V}%
-R}(w_{\rho }\left( p\left( v\right) \right) -1)\text{,}
\end{equation*}
where $n(u^{\prime })=\binom{u_{1}}{2}+\binom{u_{2}}{2}+...+\binom{u_{k}}{2}$
and $u_{i}$ is the number of vertices at distance $i$ from its roots .
\end{definition}

For the forest $F$ represented in Fig. 1, we have 
\begin{equation*}
l_{\rho _{+}}\left( F\right) =\binom{4}{2}+\binom{5}{2}+\binom{1}{2}+\left(
4-1\right) +2\left( 1-1\right) +3\left( 4-1\right) +\left( 2-1\right)
+\left( 3-1\right) +2\left( 1-1\right) =31\text{.}
\end{equation*}

\begin{theorem}
Let $n\in \mathbb{N}^{\ast }$ and $R\subset \mathbf{n}$ , $\left| R\right|
=r $, then for any ranking $\rho $ on\textbf{\ }$\mathbf{2}^{\mathbf{n}}$,
we have 
\begin{equation*}
J_{n}^{\left( r\right) }\left( q\right) =\sum\limits_{F\in \mathcal{F}%
_{n,R}}q^{l_{\rho }\left( F\right) }\text{.}
\end{equation*}
\end{theorem}

\begin{proof}
For each height $k$, between $1$ and $n-r$, and each $k$-multiplet $u=\left(
u_{1},u_{2},...,u_{k}\right) $ with $u_{i}\geq 1$ and $\left| u\right|
=u_{1}+u_{2}+...+u_{k}=n-r$, we have $\binom{n-r}{u_{1},u_{2},...u_{k}}$
possibilities of placing the $n-r$ non-root vertices on the $k$ levels, i.e.
of choosing the $k$ generations $V_{1},V_{2},...V_{k}$. So we have 
\begin{equation}
\sum\limits_{F\in \mathcal{F}_{n,R}}q^{l_{\rho }\left( F\right)
}=\sum\limits _{\substack{ k\geq 1,u_{i}\geq 1  \\ u_{1}+u_{2}+...+u_{k}=n-r
}}\binom{n-r}{u_{1},u_{2},...,u_{k}}\;q^{n(u^{\prime })}\sum_{F\in F\left(
V_{0},V_{1},...,V_{k}\right) }q^{\sum\nolimits_{v\in \mathbf{n}-R}\left(
w_{\rho }\left( p(v)\right) -1\right) }\text{.}  \tag{9.1}
\end{equation}
Here $F\left( V_{0},V_{1},...,V_{k}\right) $ is the subset of forests\ of $%
\mathcal{F}_{n,R}$, with a specified height $k$, and whose generations $%
V_{0}=R,V_{1},V_{2},...V_{k}$ are all specified. This subset of forests is
clearly in bijection with the Cartesian product $\prod\nolimits_{i=0}^{k-1}%
\mathcal{G}_{i}$ where, for $i=0$ to $k-1$, $\mathcal{G}_{i}\ $is the set of
maps from $V_{i+1}$ to $V_{i}$; the map $g_{i}\in \mathcal{G}_{i}$
associated with a forest of $F\left( V_{0},V_{1},...,V_{k}\right) $, being
defined by $g_{i}\left( v\right) =p(v)$. We then have 
\begin{equation*}
\sum\limits_{v\in \mathbf{n}-R}\left( w_{\rho }\left( p(v)\right) -1\right)
=\sum\limits_{i=0}^{k-1}\sum\limits_{v\in V_{i+1}}\left( w_{\rho }\left(
g_{i}(v)\right) -1\right) \text{,}
\end{equation*}
hence, 
\begin{equation}
\sum_{F\in F\left( V_{0},V_{1},...,V_{k}\right) }q^{\sum\limits_{v\in 
\mathbf{n}-R}\left( w_{\rho }\left( p(v)\right) -1\right)
}=\prod\limits_{i=0}^{k-1}\Theta _{i}\text{,}  \tag{9.2}
\end{equation}
with 
\begin{equation*}
\Theta _{i}=\sum\limits_{g_{i}\in \mathcal{G}_{i}}q^{\sum\nolimits_{v\in
V_{i+1}}w_{\rho }(g_{i}(v))-1}\text{.}
\end{equation*}
Let us calculate $\Theta _{i}$ for $i$ between $0$ and $k-1$. By definition
of $\rho $, $w_{\rho }=\rho \left( V_{i}\right) $ is a bijection from $V_{i}$
to $\mathbf{u}_{i}$ and $\rho \left( V_{i+1}\right) $ is a bijection from $%
V_{i+1}$to $\mathbf{u}_{i+1}$. So, for each $g_{i}\in \mathcal{G}_{i}$, $%
h_{i}=\Phi \left( g_{i}\right) =\rho \left( V_{i}\right) \circ g_{i}\circ
\left( \rho \left( V_{i+1}\right) \right) ^{-1}$ is a map from $\mathbf{u}%
_{i+1}$ to $\mathbf{u}_{i}$. By construction $\Phi $ is a bijection from $%
\mathcal{G}_{i}$ to $\mathbf{u}_{i}^{\mathbf{u}_{i+1}}$, so we can make the
bijective change of index $h_{i}=\Phi (g_{i})$ . We can therefore write 
\begin{equation*}
\Theta _{i}=\sum\nolimits_{h_{i}\in \mathbf{u}_{i}^{\mathbf{u}%
_{i+1}}}\prod\limits_{j=1}^{u_{i+1}}q^{h_{i}(j)-1}=\left(
1+q+...+q^{u_{i}-1}\right) ^{u_{i+1}}\text{,}
\end{equation*}
the second equality above, comes from classical combinatorial results (see
for example Th.A p. 127 in [5] with $M=\mathbf{u}_{i+1}$, $N=\mathbf{u}_{i}$%
, $u\left( x,y\right) =q^{y-1}$ and $\mathcal{R}=M\times N$). By
transferring these expressions into $\left( 9.1\right) $ we obtain 
\begin{equation*}
\sum\limits_{F\in \mathcal{F}_{n,R}}q^{l_{\rho }\left( F\right)
}=\sum\limits _{\substack{ k\geq 1,u_{i}\geq 1  \\ u_{1}+u_{2}+...+u_{k}=n-r
}}\binom{n-r}{u_{1},u_{2},...,u_{k}}\;q^{n(u^{\prime
})}\prod\limits_{i=0}^{k-1}\left[ u_{i}\right] _{q}^{u_{i+1}}\text{,}
\end{equation*}
which with $\left( 8.1\right) $ ends the proof.
\end{proof}

\begin{remark}
We can include the case $R=V=\mathbf{n}$ in the previous theorem. $\mathcal{F%
}_{n,V}$ reduces to the empty graph $E_{n}$ ($n$ vertices without edges) and
it suffices to set for all ranking $\rho $, $l_{\rho }\left( E_{n}\right) =0$%
. We verify that $J_{n}^{\left( n\right) }=1=\sum\nolimits_{F\in \mathcal{F}%
_{n,V}}q^{l_{\rho }(F)}=q^{l_{\rho }\left( E_{n}\right) }$.
\end{remark}

\textbf{Case of reciprocal polynomials. }By replacing the formulas $\left(
7.1bis\right) $ in the formulas of Theorem 8.2 we can obtain explicit
formulas for the $\overline{J_{n}^{\left( r\right) }}.$ For example we get
from $\left( 8.1\right) $:

\begin{corollary}
\bigskip For any couple of integers $\left( n,r\right) $ satisfying $n-1\geq
r\geq 1$ \ we have 
\begin{equation}
\overline{J_{n}^{\left( r\right) }}(q)=\sum \left[ r\right] _{q}^{u_{1}}%
\left[ u_{1}\right] _{q}^{u_{2}}...\left[ u_{k-1}\right] _{q}^{u_{k}}q^{%
\sigma (u)+r\left( n-r-u_{1}\right) }\binom{n-r}{u_{1},u_{2},...,u_{k}}\text{%
,}  \tag{9.3}
\end{equation}
where the sum is over the $k$-multiplets of strictly positive integers $%
u=\left( u_{1},u_{2},...,u_{k}\right) $ with $k\geq 1$, $\left| u\right|
=n-r $, and $\sigma \left( u\right) =\sum\nolimits_{1\leq i,j\leq k,\;j\geq
i+2}u_{i}u_{j}$ \ \ 
\end{corollary}

It is possible to define level statistics from the formula $\left(
9.3\right) $ which can still be written 
\begin{equation}
\overline{J_{n}^{\left( r\right) }}(q)=\sum \left[ u_{0}\right] _{q}^{u_{1}}%
\left[ u_{1}\right] _{q}^{u_{2}}...\left[ u_{k-1}\right] _{q}^{u_{k}}q^{%
\sigma (\widehat{u})}\binom{n-r}{u_{1},u_{2},...,u_{k}}\text{,}  \tag{9.4}
\end{equation}
where the sum is over the multiplets $\widehat{u}=\left( u_{0}=r,\text{ }%
u_{1},u_{2},...,u_{k}\right) $ such that $\left| \widehat{u}\right|
=u_{0}+u_{1}+...+u_{k}=n$ and $\sigma \left( \widehat{u}\right)
=\sum\nolimits_{0\leq i,j\leq k,\;j\geq i+2}u_{i}u_{j}$.

We let the reader verify that for any ranking $\rho $ on\textbf{\ }$\mathbf{2%
}^{\mathbf{n}}$, $\overline{l_{\rho }}$ defined below is a statistic for $%
\overline{J_{n}^{\left( r\right) }}\left( q\right) $ 
\begin{equation*}
\overline{l_{\rho }}\left( F\right) =\sigma (\widehat{u})+\sum\limits_{v\in 
\mathbf{n}-R}\left( w_{\rho }\left( p(v)\right) -1\right) \text{.}
\end{equation*}
\ 

Let us point out in conclusion that the above formulas related to the $%
J_{n}^{\left( r\right) }$, and in particular formula $\left( 6.5\right) $,
lend themseves to further combinatorial developments that will be presented
in future articles.

Moreover in article $\left[ 1\right] $, the author introduces a family of
polynomials denoted $J_{\lambda }$, indexed by the set of integer partitions
and with coefficients in $\mathbb{Z}$. These polynomials generalize the $%
J_{n}$ polynomials in a direction different from that of the $J_{n}^{\left(
r\right) }$ polynomials. It is conjectured and partially demonstrated in $%
\left[ 1\right] $ that these $J_{\lambda }$ polynomials, like the $J_{n}$,
have strictly positive and log-concave coefficients.\bigskip\ Note that
formulas $\left( 6.4\right) $ and $\left( 6.5\right) $ above are used in $%
\left[ 1\right] $.

\bigskip \textbf{References}

$\left[ 1\right] $ V. Brugidou, On a particular specialization of monomial
symmetric functions, (2023), arXiv: 2306.15300.

$\left[ 2\right] $ V. Brugidou, $q$-power symmetric functions and $q$%
-exponentiel formula, (2024), arXiv : 2401.17687.

$\left[ 3\right] $ L. Carlitz, On Abelian fields. \textit{Trans. Amer. Math.
Soc.} \textbf{35} (1933) 122-136.

$\left[ 4\right] $ L. Carlitz, $q$-Bernoulli numbers and polynomials. 
\textit{Duke Math. J.} \textbf{15} (1948) 987-1000.

$\left[ 5\right] $ L. Comtet, \textit{Advanced Combinatorics}. Springer
Netherlands, Dordrecht,1974.

[6] A. de Medicis, P. Leroux, A unified combinatorial approach for $q$- (and 
$p,q$-) Stirling numbers. \textit{J. Stat. Planning and Inference} \textbf{34%
} (1993) 89-105.

$\left[ 7\right] $ U. Duran, M. Acikgoz, Apostol type ($p,q$)-Bernoulli, ($%
p,q$)-Euler and ($p,q$)-Genocchi Polynomials and numbers. \textit{Comm. in
Math. and Appl}. \textbf{8, }vol.1\textbf{\ }(2017) 7-30.

$\left[ 8\right] $ H. W. Gould, The $q$-Stirling numbers of the first and
second kinds. \textit{Duke Math. J}. \textbf{28} (1961) 281-289.

[9] V. Kac, P. Cheung, \textit{Quantum Calculus. }Springer, New York, 2002.

$\left[ 10\right] $ R. Kenyon, M. Yin, Parking functions: From combinatorics
to probability, \textit{Methodol. Comput. Appl. Probab.} \textbf{25},
Article 32 (2023).

$\left[ 11\right] $J.P.S. Kung, C.H. Yan, Goncarov polynomials and parking
functions. \textit{J. Comb. Theory Ser. A} \textbf{102 }(1) (2003) 16-37.

$\left[ 12\right] $ I. G. Macdonald, \textit{Symmetric functions and Hall
polynomials}, second ed., Oxford University Press, New York, 1995.

$\left[ 13\right] $ C. L. Mallows and J. Riordan, The inversion enumerator
for labeled trees. \textit{Bull. Amer. Soc. }\textbf{74} (1968) 92-94.

$\left[ 14\right] $ J. W. Moon, \textit{Counting Labelled Trees}, Canadian
Math. Monographs, no. 1, Canadian Mathematical Congress, 1970.

$\left[ 15\right] $ B. E. Sagan, J.P. Joshua, $q$-Stirling numbers in type B
(2022), preprint, arXiv: 2205.14078

$\left[ 16\right] $ R.P. Stanley, Hyperplane arrangements, parking
functions, and tree inversions. In\textit{\ Mathematical Essays in Honor of
Gian-carlo Rota} \textit{(Cambridge, MA, 1996)}, vol. 161 of Progr. Math.,
Birkh\"{a}user, Boston, MA (1998) 359-375.

$\left[ 17\right] $ R. P. Stanley, \textit{Enumerative Combinatorics,} 
\textit{Vol. 2}. Cambridge University Press, Cambridge, UK, 1999.

$\left[ 18\right] $ R. P. Stanley, Some aspects of (r,k)-parking functions. 
\textit{J. Comb. Theory Ser. A} \textbf{159} (2018) 54-78.

[19] M. Wachs and D. White, p,q-Stirling numbers and set partition
statistics. \textit{J. Comb. Theory Ser. A} \textbf{56} (1991) 27-46.

$\left[ 20\right] $ C. H. Yan, Generalized parking functions, tree
inversions, and multicolored graphs. \textit{Adv. in Appl. Math}. \textbf{27 
}(2-3) (2001) 641-670.

[21] C.H. Yan, Parking functions, in M. Bona (ed.), \textit{Handbook of
Enumerative Combinatorics}, CRC Press, Boca Raton, FL (2015) pp. 835-893.

\end{document}